\documentclass[11pt]{article}%
\usepackage{amsmath}
\usepackage{amsfonts}
\usepackage{amssymb}
\usepackage{graphicx}
\usepackage{theorem}
\usepackage{color}%
\providecommand{\U}[1]{\protect\rule{.1in}{.1in}}
\newtheorem{thm}{Theorem}[section]

\newtheorem{rmk}[thm]{Remark}

{\theorembodyfont{\upshape}
\newtheorem{examp}[thm]{Example}
}
\numberwithin{equation}{section} 
\begin{document}

\title{Title: }

\begin{center}
\textbf{SEQUENTIAL\ CONE--COMPACTNESS\ DOES\ NOT\ IMPLY\ CONE--COMPACTNESS}

\bigskip

by

\bigskip

Marius\ DUREA\footnote{{\small Faculty of Mathematics, \textquotedblright
Alexandru Ioan Cuza\textquotedblright\ University, 700506--Ia\c{s}i, Romania
and \textquotedblright Octav Mayer\textquotedblright\ Institute of
Mathematics, Ia\c{s}i Branch of Romanian Academy, 700505--Ia\c{s}i, Romania{;
e-mail: \texttt{durea@uaic.ro}}}} and Elena-Andreea
FLOREA\footnote{{\small Faculty of Mathematics, \textquotedblright Alexandru
Ioan Cuza\textquotedblright\ University, 700506--Ia\c{s}i, Romania;{ e-mail:
\texttt{andreea.florea@uaic.ro}}}}

\bigskip

\bigskip
\end{center}

\noindent{\small {\textbf{Abstract:}} We address a problem posed in [1] by
demonstrating through an example that, in the absence of separability, the
property of sequential cone compactness does not generally imply cone
compactness.}

\bigskip

\noindent{\small {\textbf{Keywords:}} cone compactness{ $\cdot$ sequential
cone compactness $\cdot$ separability}}

\bigskip

\noindent{\small {\textbf{Mathematics Subject Classification (2020): }}54A20}

\begin{center}

\end{center}

\section{Preliminaries}

Let $X$ be a normed vector space, $C\subset X$ be a closed convex cone, and
$A\subset X$ be a nonempty set. According to \cite{Luc} the set $A$ is called
$C-$compact (or compact with respect to the cone $C$) if from any cover of $A$
with the sets of the form $U+C$, where $U$ is open, one can extract a finite
subcover of it. In \cite{DF} the following concept was introduced and studied:
the set $A$ is called $C-$sequentially compact if for any sequence $\left(
a_{n}\right)  \subset A$ there is a sequence $\left(  c_{n}\right)  \subset C$
such that the sequence $\left(  a_{n}-c_{n}\right)  $ has a convergent
subsequence towards an element of $A$. These two notions are also known under
the names of cone compactness and sequential cone compactness, respectively.

It was shown in \cite[Theorem 2.7]{DF} that if $A$ is $C-$compact, then $A$ is
$C-$sequentially compact. For the converse, according to \cite[Theorem
2.11]{DF}, if $A$ is $C-$sequentially compact and separable, then $A$ is
$C-$compact. The question of whether the separability assumption can be
omitted was left unanswered in that paper (see \cite[Remark 2.12]{DF}). The
aim of this technical note is to demonstrate that separability is indeed
essential. To this end, we provide an example in a nonseparable normed vector
space of a set that is $C-$sequentially compact with respect to a given cone
$C$, yet not $C-$compact.

Even without a clear statement outlining the full relationship between these
concepts in a general normed vector space, sequential cone compactness has
already been successfully employed as a technically more manageable substitute
for cone compactness (see \cite{H} and the references therein). Therefore, by
demonstrating that, beyond the realm of separability, these two notions are
indeed distinct---showing that $C-$sequential compactness is strictly weaker
than $C-$compactness---this example helps clarify several results in the literature.

\section{The main result}

The next example shows that sequential cone compactness does not imply cone compactness.

\begin{examp}
\label{ex}Let $X:=B\left(  \mathbb{R}\right)  $ be the Banach space of bounded
functions from $\mathbb{R}$ to $\mathbb{R}$ with the supremum norm. Take
$C:=B\left(  \mathbb{R}\right)  _{+},$ that is the closed, convex, pointed
cone of nonnegative functions in $B\left(  \mathbb{R}\right)  .$ Consider now
the set $A\subset B\left(  \mathbb{R}\right)  $ which consists of those
functions $f$ having the property that there is $X_{f}\subset\mathbb{R}$ with
$\operatorname*{card}X_{f}=\aleph_{0}$ such that $f\left(  x\right)  =-1,$ if
$x\in X_{f},$ and $f\left(  x\right)  =0,$ if $x\in\mathbb{R}\setminus X_{f}.$
There are $2^{\aleph_{0}}$ such functions and the norm of the difference of
any two of them is $1.$ In particular, $A$ is not separable (as a metric
space). Clearly,
\[
A\subset\bigcup_{f\in A}U\left(  f,1\right)  +C,
\]
where $U\left(  f,1\right)  $ is the open ball centered at $f$ of radius $1.$
Suppose, by way of contradiction, that there exists a finite subcover of this
cover of $A.$ Therefore, there is a nonzero natural number $k$ and $\left(
f_{i}\right)  _{i\in\overline{1,k}}\subset A$ such that%
\[
A\subset\bigcup_{_{i\in\overline{1,k}}}U\left(  f_{i},1\right)  +C.
\]
Consider the countable set
\[
Y:=\bigcup_{_{i\in\overline{1,k}}}X_{f_{i}},
\]
$\overline{x}\in\mathbb{R}\setminus Y,$ and $Z=Y\cup\left\{  \overline
{x}\right\}  .$ Denote by $f$ the function from $A$ defined by the equality
$Z=X_{f}$ (that is $f\left(  x\right)  =-1,$ if $x\in Z$, and $f\left(
x\right)  =0,$ otherwise). Consequently, one should have
\[
f\in\bigcup_{_{i\in\overline{1,k}}}U\left(  f_{i},1\right)  +C,
\]
meaning that there exist $i\in\overline{1,k}$ and $h\in C$ such that
$\left\Vert f-f_{i}-h\right\Vert <1.$ But $f\left(  \overline{x}\right)  =-1$
and $f_{i}\left(  \overline{x}\right)  =0,$ so, because $h\left(  \overline
{x}\right)  \geq0,$
\[
1\leq\left\vert f\left(  \overline{x}\right)  -f_{i}\left(  \overline
{x}\right)  -h\left(  \overline{x}\right)  \right\vert \leq\left\Vert
f-f_{i}-h\right\Vert <1,
\]
which is a contradiction. We deduce that one cannot extract a finite subcover
from the above conic cover of $A,$ whence $A$ is not $C-$compact. (Actually,
the same argument shows that one cannot extract even a countable subcover of
the given cover.)

Consider now a sequence $\left(  f_{n}\right)  _{n\geq1}\subset A.$ Take the
countable set
\[
X:=\bigcup_{n\geq1}X_{f_{n}},
\]
and define the function $g\in A$ by the equality $X=X_{g}\ $(that is,
$g\left(  x\right)  =-1,$ for all $x\in X$ and $g\left(  x\right)  =0,$
otherwise). Now consider the sequence $\left(  g_{n}\right)  _{n\geq1}$ of
functions defined by $g_{n}=f_{n}-g,$ for all $n\geq1.$ It is easy to see that
$\left(  g_{n}\right)  _{n\geq1}\subset C$ and for all $n\geq1,$ $f_{n}%
-g_{n}=g\in A.$ Therefore, $A$ is $C-$sequentially compact.
\end{examp}

\begin{rmk}
As one can observe, the set $A$ in Example \ref{ex} possesses a much stronger
property than what is required by the definition of $C-$sequential
compactness. Specifically, for every sequence $\left(  f_{n}\right)  _{n\geq
1}\subset A$ there exists a sequence $\left(  g_{n}\right)  _{n\geq1}\subset
C$ such that $\left(  f_{n}-g_{n}\right)  _{n\geq1}$ is stationary at some
point in $A.$ However, $A$ is not $C-$compact.
\end{rmk}

\bigskip

\noindent\textbf{Data availability.} This manuscript has no associated data.

\noindent\textbf{Disclosure statement}. No potential conflict of interest was
reported by the authors.

\end{document}